\documentstyle[12pt]{article} 
\normalbaselines

\begin{document}

\newtheorem{lema}{Lemma}[section]
\newtheorem{defi}[lema]{Definition}  
\newtheorem{rema}[lema]{Remark}  
\newtheorem{cont}[lema]{Counterexample}
\newtheorem{prop}[lema]{Proposition}
\newtheorem{teor}[lema]{Theorem}
\newtheorem{coro}[lema]{Corollary}
\newtheorem{ejem}[lema]{Example}  

\font\ddpp=msbm10 scaled \magstep 1 
\def\R{\hbox{\ddpp R}}    
\def\C{\hbox{\ddpp C}}     
\def\L{\hbox{\ddpp L}}    
\def\N{\hbox{\ddpp N}}    
\def\Z{\hbox{\ddpp Z}}    
\def\H{\hbox{\ddpp H}}    
\def\S{\hbox{\ddpp S}}    

\renewcommand{\theequation}{\thesection.\arabic{equation}} 

\begin{center}
{\large \bf SOME SEMI-RIEMANNIAN VOLUME COMPARISON THEOREMS}\\
Paul E. Ehrlich and Miguel S\'anchez

\end{center}

\begin{abstract}
\noindent  Lorentzian versions of classical Riemannian volume comparison theorems by Günther, Bishop and Bishop-Gromov, are stated for suitable natural subsets
of general semi-Riemannian manifolds. The problem is more subtle in the 
Bishop-Gromov case, which is extensively discussed. For the general
semi-Riemannian case, a local version of the Günther and Bishop theorems is
given and applied. 
\end{abstract} 

\vspace{1mm}

\bibliographystyle{unsrt}

\section{Introduction}

A substantial topic in Riemannian Geometry is the comparison between the volumes of subregions of (complete) 
Riemannian manifolds under various curvature hypotheses. Among the classical results are: 
(1) Günther's theorem \cite{Gu},  which imposes an upper bound $c$ for the sectional curvature of a 
Riemannian manifold, and obtains a lower bound for the volume of geodesic balls in the manifold in terms of 
the volume of balls with equal radii in the model space of curvature $c$, (2) Bishop's theorem \cite{B}, 
\cite{BC}, which imposes a lower bound for the Ricci curvature (in the spirit of Myers' theorem) and 
obtains a reversed inequality for the volumes of corresponding metric balls, and (3) the Bishop-Gromov theorem
\cite{G} which, under the hypotheses in Bishop's theorem, obtains the non-increasing monoticity of 
the ratio between the volumes of metric balls in the manifold and the model space. A detailed account of these 
results may be found, for example, in \cite{Ch} or \cite{GHL}. 

\vspace{3mm}
\noindent Recently, some comparison results between volumes of semi-Riemannian manifolds have been obtained 
(see \cite{AH}, \cite{EJK}  and references therein). There are several difficulties to extend the Riemannian 
Bishop-type theorems to indefinite  semi-Riemannian manifolds. For example, by results of Kulkarni \cite{K} and  
Dajczer and Nomizu \cite{DN}  if an (upper or lower) bound is imposed on the sectional 
(or, in an appropriate sense, Ricci) curvature of an indefinite manifold of dimension $\geq 3$, then the curvature 
must be constant. Thus, it is natural to restrict the curvature inequality to planes with a definite 
causal character \cite{EJK} or to make an inequality directly on the curvature tensor $<R(x,y)y,x>$ itself, 
cf. \cite{AH}. 
The following further two difficulties also influence our approach to volume comparison:

(i) For indefinite manifolds, the metric balls make sense only in the Lorentzian case, where the 
``Lorentzian distance" (or ``time-separation") may be defined. Nevertheless, even in this case these balls 
may behave very differently than Riemannian metric balls. For example, the ``inner" metric balls $B^+(p,\epsilon) =
\{ q \in I^+(p); d(p,q) < \epsilon \}$ need not be open in general space-times since the space-time distance 
function may fail to be upper semicontinuous, cf. \cite{BEE}, p. 142. Further, the analogues of the 
closed Riemannian balls for space-times may have infinite volume, and are usually neither compact nor contained in any normal neighborhood. This suggests studying the relationship between the 
volumes of transplanted subsets included in normal neighborhoods rather than the metric balls themselves
(see Section 2 for detailed definitions). 

(ii) It is possible to find indefinite manifolds (notably, any pair of model spaces with different curvature, 
see Corollary \ref{c}) and two normal neighborhoods $W, W'$ of a point $p$ in the first model space such 
that the corresponding transplanted neighborhoods $\hat{W}, \hat{W}'$ of $\hat{p}$ in the second manifold satisfy 
vol$(W) <$ vol$(\hat{W})$ and vol$(W') >$ vol$(\hat{W}')$. These inequalities are a consequence of the different 
behavior of timelike and spacelike directions under bounds on the sectional or Ricci curvatures.

\vspace{3mm}
\noindent These above considerations lead us to deal with special neighborhoods which make possible volume 
comparison. We will chose two such types of neighborhoods: (a) in the Lorentzian case, those we will 
call {\it SCLV} at a point $p$, which are subsets of a normal neighborhood of $p$ covered by the timelike 
geodesics emanating from $p$, and (b) in the general semi-Riemannian case, ``small" neighborhoods 
around either a spacelike or a timelike direction.

SCLV subsets are defined in Section 2. For this definition, we do not need any global causal assumption 
on the manifold; nevertheless, if the Lorentzian manifold is globally hyperbolic, then 
some canonical restrictions can be done, and it seems natural to consider SCLV subsets which are contained in certain subsets ${\cal C}^+(p)$ as is done in \cite{EJK} (cf. the discussion at the beginning of Section 2). In Section 2, Theorems \ref{t1}, \ref{t2}, \ref{t3} are also stated; they are the Lorentzian versions of, respectively, the Günther-Bishop, the Bishop 
and the Bishop-Gromov theorems for SCLV subsets. Their proofs are given in Section 3 and may be given similarly
to the Riemannian comparison theory, since the metric tensor restricted to the orthogonal complement of a 
timelike (radial) geodesic is positive definite. In considering the general notion of SCLV subsets in the context of
volume ratio comparison analogous to that of the Bishop-Gromov
Theorem, in Counterexample \ref{3.3} a special property of the Riemmanian
metric balls or of the distance wedges of \cite{EJK} emerges, which is
partly responsible for this monotonicity: these sets
traditionally studied in the previous references have the property
that their cut functions (in the sense of (\ref{e1}) and Condition (B) of
Theorem \ref{t3}) are constant. On the other hand, if we restrict our attention
to volume comparison with the Lorentz-Minkowski spacetime, which naturally arises in the case often considered
in General Relativity that (M,g) is globally hyperbolic and satisfies the Timelike Convergence Condition
Ric$(v,v) \geq 0$ for all timelike $v$, then a volume comparison for general subsets $Z$ of $J^+(p)$ is obtained
in Corollary \ref{c-t3}. 

In the Riemannian case, the existence of ``many" manifolds admitting any (non-trivial) type of bounds for 
the Ricci or sectional curvatures is elementary. For example, if at a point $p$ the Ricci curvature is 
greater than $c$, then small variations of the metric (and its two first derivatives) around $p$ yield 
new metrics satisfying the same inequality. The question is not so transparent in the Lorentzian case, 
because of the previously recalled results \cite{DN}, \cite{K} that certain curvature inequalities imply
constant curvature. Moreover, the lack of compactness of the unit tangent vectors at a point suggests that
explicit calculations should be made. Thus in Section 4, some examples of Lorentzian manifolds, which satisfy
the curvature bounds of our volume comparison theorems, are given. The first family of examples consists 
of warped products with a negative interval as basis which, following \cite{ARS-grg}, are called Generalized 
Robertson-Walker. The second family is obtained by changing conformally a given Lorentzian metric. The 
conformal factors used are ``radial" and, so, the obtained curvature bounds apply just on radial planes 
or directions, which is the minimum hypothesis we need. The comparisons in Theorems \ref{t1} and \ref{t2} 
are carried out with model spaces just for simplicity, and could be stated between arbitrary Lorentzian 
manifolds under the corresponding inequalities for sectional or Ricci curvatures; conformal changes yield 
examples of this situation also.

In Section 5 semi-Riemannian manifolds of arbitrary index are considered, and we obtain local versions of 
the Günther and Bishop theorems for the so called {\it SCV } subsets, in Theorem \ref{t}. As an application 
to the Riemannian case, we give a local version of the Günther-Bishop Theorem, where the bound on the 
sectional curvature is replaced by one on the Ricci curvature, in Corollary \ref{c-tf}.

\section{The Lorentzian case: set-up and statement of the results}

Let $(M,g)$ be a time-oriented Lorentzian $n-$manifold, $n\geq 2$, and $Q_{c}$ the Lorentzian model space 
of constant curvature $c$ for some fixed $c \in \R{}$;  choose $p\in M$ and $p_0 \in Q_c $. 
A neighborhood $V$ of $p$ is normal if the exponential map $\exp_p$ is a diffeomorphism from 
$\overline{V} = \exp_p^{-1}(V)$ to $V$, and $\overline{V}$ is starshaped from the origin $0_p$. 
Take a $\overline{U} \subset T_p M$ satisfying: 
\begin{enumerate}

\item $\overline{U}$ is an open subset of the causal future of the origin, $\overline{U} \subset J^+(0_p)$.

\item $\overline{U}$ is starshaped from the origin (if $u\in \overline{U}$ then $tu \in \overline{U}, 
\forall t\in (0,1)$),  the exponential map at $p$, $\exp_p$, is defined on $\overline{U}$, and the 
restriction of $\exp_p$ to $\overline{U}$ is a diffeomorphism onto its image $U= \exp_p(\overline{U})$. 

\item The closure of $\overline{U}$, $cl(\overline{U})$, is compact.\footnote{This assumption is 
introduced just to avoid infinite integrals, which will be necessary for the discussion of the 
equality in the three theorems stated in this section (and will be especially relevant for 
the third one). Nevertheless, it can be weakened: with the notation to be introduced, it is 
sufficient to assume that the integral $\int _{\overline{U}} \mid s_c^{n-1} \mid$ is finite. 
In this case, note that the cut functions (2.1) may achieve the value $\infty$.}

\end{enumerate}
The subset $U=\exp_p \overline{U}$ satisfying these three conditions will be called {\it standard 
for comparison of Lorentzian volumes} (SCLV) {\it at $p$} (but $p$ is never in the interior of $U$). 
Note that if item 1 is suppressed then this definition makes sense in the general semi-Riemannian case; 
we will call {\it standard for comparison of volumes} (SCV) {\it at $p$} the subsets satisfying 2 and 3. 
In fact, all our results in this section can be extended to the Riemannian case (but not to other indices) 
just by considering neighborhoods SCV at $p$. 

{\it Example.} For each $p\in M$, consider the subset 
\[
{\cal C}^+(p)= \{ \gamma_\xi(t): g(\xi,\xi)=-1, \xi\in T_pM, \xi\in I^+(0), 0<t<{\rm cut}_p(\xi)\},
\]
where cut$_p(\xi)$ is the value of the cut function at $p$ in the direction $\xi$, and $\gamma_\xi$ is 
the inextendible geodesic with initial velocity $\xi$; that is, the Lorentzian distance between $p$ 
and $\gamma_{\xi}(t)$ is equal to $t$ for all $t\leq c_{\xi}(t)$ and strictly bigger for $t>c_{\xi}(t)$ 
(the last inequality, when $c_{\xi}(t)<\infty$ and $\gamma_{\xi}$ is defined beyond). Clearly, 
${\cal C}^+(p)$ satisfies conditions 1 and 2. When $(M,g)$ is  globally hyperbolic, by a well-known 
Avez-Seifert result (\cite{Av}, \cite{Se}),  each pair of causally related points may be joined by a 
maximizing causal geodesic. Thus, the closure of ${\cal C}^+(p)$ is $J^+(p)$, and these two subsets 
are equal, up to a zero-measure subset. In this case, it seems natural to consider subsets 
SCLV at $p$ included in ${\cal C}^+(p)$ (see \cite{EJK}), even though we do not impose this 
assumption. On the other hand, when $(M,g)$ is not globally hyperbolic ${\cal C}^+(p)$ may be empty as for so-called totally vicious spacetimes.

The {\it cut function} $c_U$ of the SCLV subset $U$ is:

\begin{equation}
\label{e1}
c_U: \overline{U}[1] \rightarrow \R{}, \quad c_U(\xi) = {\rm Sup}\{t\in(0,\infty): t\xi \in \overline{U}\},
\end{equation}
where $\overline{U}[1]$ is the set of timelike unit vectors in $\overline{U}$.

Let $i: T_pM \rightarrow T_{p_0} Q_c $ be a linear isometry, and define the usual transplantation 
$F:U\rightarrow Q_c, F=\exp_{p_0} \circ i \circ  (\exp_p \mid_U )^{-1}$. Put $\overline{U}_0 = i(\overline{U})$ 
and $U_0 = \exp_{p_0}(\overline{U}_0) (=F(U))$ (note that $\overline{U}_0$ is open, but $U_0$ not necessarily), Fig. 1.

\vspace{2.4cm}
\begin{center}
 Figure 1\\
Transplantation $F$ of a subset SCLV at $p$.\\
{\bf (See last page)}
\end{center}
\vspace{4.7cm}

A tangent vector $v$ to $U$ will be called {\it radial} if it can be written as 
$\frac{d}{dt}\mid_{t_0} \exp_p(tv_p)$ for some vector $v_p$ tangent at $p$. 
A tangent plane $\pi$ to $U$ will be {\it radially timelike} if it contains a timelike radial vector. 
The sectional curvature of $\pi$ will be denoted by $K(\pi )$. 


\begin{teor}
\label{t1}
Let $(M,g)$ be a Lorentzian manifold, let $U$ be a subset SCLV at $p \in M$, and assume that 
the following two conditions hold:

{\rm (1)} $K(\pi) \geq c$ for all radially timelike planes $\pi$.

{\rm (2)} $\exp_{p_0}: \overline{U}_0\rightarrow U_0$ is a diffeomorphism (thus, $U_0$ is assumed to be open in $Q_c$) .

\noindent Then the relation between the volumes of $U$ and $U_0$ is:
\begin{equation}
\label{e-t1}
{\rm vol}(U) \geq {\rm vol}(U_0)
\end{equation}
and the equality holds if and only if $F:U\rightarrow U_0$ is an isometry.
\end{teor}

{\it Remark.} If $c=0$, the condition (2) is trivially satisfied. If $c>0$ it can be dropped also, because 
the future timelike cut locus is empty for $Q_c$. If $c<0$, conjugate points appear along timelike 
geodesics of length greater or equal to $\pi / \sqrt{-c}$. Thus, condition (2) can be replaced by 

(2') {\it If $c<0$ then the diameter of $\overline{U}_0 $ (or the diameter of $U_0$) for the 
Lorentzian distance satisfies:} ${\rm diam}(\overline{U}_0) < \pi /\sqrt{-c}.$

\begin{teor}
\label{t2}
Let $(M,g)$ be a Lorentzian manifold, and assume that 
\[
{\rm Ric}(v,v) \geq (n-1)c \cdot g(v,v)
\]
for all timelike and radial vectors $v$ tangent to a subset $U$, SCLV at $p \in M$. Then:
\begin{equation}
\label{e-t2}
{\rm vol}(U) \leq {\rm vol}(U_0)
\end{equation}
and the equality holds if and only if $F:U\rightarrow U_0$ is an isometry.
\end{teor}

{\it Remark.} The assumption in Theorem \ref{t2} is analogous to the assumption (1) in 
Theorem \ref{t1} but the assumption on the sectional curvature is weakened to an assumption 
on the Ricci curvature, in the spirit of the Bishop Theorem. On the other hand, the assumption 
(2) in Theorem \ref{t1} is now automatically satisfied because, under the assumption in Theorem \ref{t2},
the timelike conjugate points in $U \subset M$ appear before the timelike conjugate points 
in $U_0 \subset Q_c$ (recall we assume as a general hypothesis that $\exp_p : \overline{U}\rightarrow U$ 
is a diffeomorphism and, so, no timelike cut points appear in $U_0$).

\vspace{3mm}

\noindent For a theorem analogous to the Bishop-Gromov Theorem, we will introduce the 
following notation. For each $r>0$ put $\overline{U}^r = r \overline{U} (=\{ ru: u\in \overline{U}\}), 
\overline{U}^r_0 = r \overline{U}_0 , U^r = \exp_p(\overline{U}^r), U^r_0=\exp_{p_0}(\overline{U}^r_0)$. 
Take $I=(0,b), b\in(0,\infty]$ such that the $\overline{U}^r$ are also subsets SCLV  at $p, \forall r\in I$. 
(Note that $b$ can be always taken at least equal to 1.)
  
\begin{teor}
\label{t3}
Let $(M,g)$ be a Lorentzian manifold, and assume as in Theorem \ref{t2} that 
\[
{\rm Ric}(v,v) \geq (n-1)c \cdot g(v,v)
\]
for all timelike and radial vectors $v$ tangent to a subset $U,$ SCLV at $p\in M.$ The function 
$V: I\rightarrow \R{}^+$, $V(r) = {\rm vol}(U^r) / {\rm vol}(U_0^r)$ is non-increasing, if 
one of the following two conditions hold:

(A) c=0

(B) The cut function $c_U$ of $U$ (defined in formula {\rm (\ref{e1})}) is constant. 

\noindent Moreover, in this case if there exist $r,R\in I, r<R$ such that $V(r)=V(R)$ then $U^R$ is 
isometric to $U_0^R =F(U^R)$.

\end{teor}

{\it Remark.} The necessity of one of the conditions {\em (A) {\em or} (B)} will be discussed after 
the proof in the next Section. On the other hand, recall that both conditions make sense in the 
Riemannian case for any neighborhood SCV at $p$. In fact, condition {\it (B)} says that 
$\overline{U}= R_0 \cdot \overline{U}[1]$, where $c_U \equiv R_0$, and $\overline{U}[1]$  is now 
the set of unit (spacelike) vectors in $\overline{U}$. When $\overline{U}[1]$ is the whole unit 
sphere in $T_pM$ then $U$ is a ball of radius $R_0$, and our theorem is just the (Riemannian) Bishop-Gromov 
theorem. Thus, condition {\it (A)} yields an extension of this classical theorem for comparison with 
Euclidean space. Condition {\em (B)} also yields an extension of \cite[Theorem 4.4]{EJK}.

\vspace{3mm} 
\noindent  Condition {\em (A)} is especially interesting from the point of view of the General Relativity. 
In fact, this assumption is just ${\rm Ric}(v,v)\geq 0$ for all timelike radial vectors, and the 
well-known relativistic Timelike Convergence Condition (which means that the gravity, on average, attracts) 
imposes that the Ricci curvature be $\geq 0$ on all timelike vectors. Moreover, we can identify $Q_c$ for 
$c=0$ with $T_pM$, and $\overline{U} \equiv U_0$.  When $(M,g)$ is globally hyperbolic, we can canonically 
assume that the exponential map is a diffeomorphism $\exp_p: \overline{{\cal C}}^+(p) \rightarrow {\cal C}^+(p)$. In the following consequence of the proof of Theorem \ref{t2} and Theorem \ref{t3}, $Z^*$ denotes the 
smallest starshaped subset of ${\cal C}^+$ containing $Z\cap {\cal C}^+$,  for each $Z\subset J^+(p)$.

\begin{coro}
\label{c-t3}
Let $(M,g)$ be a globally hyperbolic Lorentzian manifold
satisfying the Timelike Convergence Condition, and $p\in M$.

For any subset $Z \subseteq J^+(p)$ in $M$ we have:
\[
{\rm vol} (Z) \leq {\rm vol} (Z_0), \quad {\rm where} \quad Z_0=\exp_p^{-1}(Z) \subset \overline{{\cal C}}^+(p) ,  
\]
and, if ${\rm vol}(Z) < \infty$, the equality holds if and only if $Z^*$ is flat.

For any subset $U,$ SCLV at $p\in M$, the function 
\[
V: I\rightarrow \R{}^+, V(r)=
{\rm vol}(U^r) / {\rm vol}(U_0^r)
\]
is non-increasing. Moreover, 
there exist $r,R\in I, r<R$ such that $V(r)=V(R)$ if and only if  $U^R$ is flat.

\end{coro}

\section{Proofs of the theorems}

Our convention for the curvature tensor will be:
$R(X,Y)= [\nabla_X , \nabla_Y] - \nabla_{[X,Y]}.$ The following Lorentzian version of Rauch's Comparison Theorem will be needed:

\begin{teor}
\label{Rauch}
Let $(M_1,g_1), (M_2,g_2)$ be two time-oriented Lorentzian manifolds with dim$M_1 \leq $ dim$M_2$. Let $\gamma_i:[0,b] \rightarrow M_i$ be future directed timelike unit geodesic segments, $i=1,2$. Put $p_i = \gamma_i(0)$, fix a linear isometry $j:T_{p_1}M_1 \rightarrow T_{p_2}M_2 $ and assume that for each $v_1 \in T_{p_1}M_1$ orthogonal to $\gamma_1 '(0)$ and $v_2= j(v_1)$, the following relation between sectional curvatures at $M_1, M_2$ holds:
\begin{equation}
\label{e-Rauch}
K_1 (\pi_1(t) ) \geq K_2 (\pi _2(t) ), \quad \forall t\in (0,b)
\end{equation}
where each plane $\pi_i(t)$ is the parallel transport along $\gamma_i (t)$ of Span$\{ v_i, \gamma'_i(0)\}$. Let $J_i$ be Jacobi fields on $\gamma_i$ orthogonal to $\gamma'_i$ and such that $J_i(0)=0, j(J'_1(0))=J'_2(0)$, and put $\phi(t)= g_1(J_1,J_1)(t)/g_2(J_2,J_2)(t)$. 

If $\gamma_1 $ has no conjugate points to $t=0$ on $(0,b)$, then $\lim_{t\rightarrow 0} \phi (t) =1$ and
$\phi'(t)\geq 0$.
Moreover, if $\phi (t_0) = 1$ for some $t_0\in (0,b]$ then the equality holds in (\ref{e-Rauch}), $\forall t \in [0,t_0]$ choosing $v_1 =J_1'(0)$. 
\end{teor}

This version is straightforward from the techniques in \cite[p. 149]{BEE} and references therein.
 
\vspace{3mm}
\noindent {\bf Proof of Theorem \ref{t1}.} For each timelike unit $\xi \in \overline{U}[1]$ consider the 
inextendible geodesic $\gamma_\xi : I_{\xi}\rightarrow U$ with initial velocity $\xi$ and $I_{\xi}=[0,b_\xi)$. 
Note that there are no conjugate points in $I_{\xi}$ and, for each $t\in I_{\xi}$ define the endomorphism of 
the orthogonal subspace to $\xi$:
\[
{\cal A}_{\xi}(t): \xi ^{\perp}\rightarrow \xi ^{\perp}, v\rightarrow (\tau_0^t)^{-1}(J^{\xi}_v(t))
\]    

\noindent where $\tau_{0}^{t}$ is the parallel transport from $0$ to $t$ along $\gamma_{\xi}$ and $J^{\xi}_{v}$ 
denotes the unique Jacobi field along $\xi$ with zero initial value, $J^{\xi}_{v}(0)=0$, and covariant 
derivative $ J^{\xi '}_{v}(0)= v. $ As in the Riemannian case, ${\cal A}_{\xi}$ satisfies:

(1) ${\cal A}_{\xi}(0)=0$, and it is differentiable in $t$.

(2) ${\cal A}_{\xi}^{(k)}(t)(v)=(\tau_0^t )^{-1}(J^{\xi}_{v})^{(k)}(t)$, where $^{(k)}$ denotes the 
$k$-th covariant derivative. In particular, ${\cal A}_{\xi}'(0)=I$ (identity).

(3) ${\cal A}_{\xi}''(t) + R_{\xi}(t){\cal A}_{\xi}(t) = 0$, where $R_{\xi}(t)(v) =
(\tau_0^t)^{-1}[R(\tau_0^t(v),\gamma'_{\xi}(t))\gamma'_{\xi}(t)]$
which is also an endomorphism of $\xi^{\perp}$ varying differentiably with $t$.

\vspace{3mm}
\noindent Consider the functions:
\[
s_c(t)=\left\{ 
\begin{array}{cl}
\frac{\sin(\sqrt{-c} t)}{\sqrt{-c}} & \quad {\rm if} \quad c<0 \\
t				    & \quad {\rm if} \quad c=0 \\
\frac{\sinh(\sqrt{c} t)}{\sqrt{c}} & \quad {\rm if} \quad c>0 
\end{array} \right. \quad \quad 
c_c(t)=\left\{ 
\begin{array}{cl}
\cos(\sqrt{-c} t) & \quad {\rm if} \quad c<0 \\
1                 & \quad {\rm if} \quad c=0 \\
\cosh(\sqrt{c} t) & \quad {\rm if} \quad c>0 
\end{array} \right.
\] 

\noindent which satisfy: $s'_c(t)=c_c(t), c'_c (t)=c s_c(t) , c^2_c-cs_c^2=1.$ Note that if $(M,g)$ has constant curvature $c$ then ${\cal A}_{\xi}(t)=s_c(t) I$ (if $V$ is a parallel vector field on $\xi$ then $J(t)=s_c(t)V$ is a Jacobi field, with $J''(t)=cJ, R(J,\gamma'_{\xi})\gamma'_{\xi}=-cJ)$.

Now, let us  study ${\rm vol}(U)$. First, note that as $U$ is a SCLV subset: 
\[
{\rm vol}(U) = \int_{\overline{U}} {\rm Jac}(exp_p) = \int_{\overline{U}[1]} d\xi \int_0^{c_u(\xi)} dt \, \, {\rm det}({\cal A}_{\xi}(t)),
\]

\noindent where for the last equality, we have used: (i) in general, if no conjugate point has appeared, ${\rm Jac}(exp_p)(t\xi ) = {\rm det}({\cal A}_{\xi}(t))$, (ii) the null cone has 0-measure, and (iii) in a SCLV subset, polar integration can be applied.

Thus, we have just to prove that, under our hypotheses, 
\begin{equation}
\label{edes}
{\det} {\cal A}_{\xi}(t) \geq s_c(t)^{n-1}, \quad \forall t \in (0,c_U(\xi)).
\end{equation}
to obtain the inequality between the volumes (\ref{e-t1}). Moreover, 
we will prove that the equality in (\ref{edes}) holds at $t_0 \in (0,c_U(\xi))$ if and only if the sectional curvature of the (timelike) radial planes along $\gamma_{\xi}\mid_{(0,t_0)}$ is equal to $c$; thus, 
the discussion of the equality of the volumes is a straightforward consequence of the semi-Riemannian version of Cartan's theorem (see, for example, \cite[p. 222-3]{O} taking into account that the the version of Theorem 14 there can be extended as in the Riemannian case, in the spirit of \cite[Ch. 7]{Sp}; see also \cite[1.7.18]{Wo} for affine Cartan's theorem). Summing up, the following result will finish the proof.

\begin{lema}
\label{l1}
Put $\psi (t) \equiv {\det} {\cal A}_{\xi}(t) / s_c(t)^{n-1}$. Under the hypotheses of Theorem \ref{t1}:
\[
\psi(t)  \geq 1, \quad \forall t \in (0,b_{\xi}).
\]
Moreover, if the equality holds at $t_0 \in (0,b_{\xi})$ then:
\[
 {\cal A}_{\xi}(t) = s_c(t) I, \quad {\rm {\it and}} \quad R_{\xi}(t)=-cI
\quad \forall t \in [0,t_0].
\]
\end{lema}

{\it Proof.} Clearly, by L'Hopital's rule $
\lim_{t\rightarrow 0} \psi (t)= \lim_{t\rightarrow 0} 
{\rm det} {\cal A}'_{\xi}(t) =1. $
Thus, it suffices to prove $\psi'(t) \geq 0$ on $(0,b_{\xi})$, or, equivalently:
\begin{equation}
\label{e-l1-1}
\frac{({\rm det} {\cal A}_{\xi})'(t)}{{\rm det}{\cal A}_{\xi}(t)} \geq (n-1) \frac{s'_c(t)}{s_c(t)}. 
\end{equation}
Now, put ${\cal B}(t) \equiv {\cal B}_{\xi}(t) = \hat{{\cal A}}_{\xi}(t) \circ {\cal A}_{\xi}(t)$, where the symbol $\hat{}$ denotes adjoint. Then ${\cal B}(t)$ is invertible close to 0, and (\ref{e-l1-1}) is equivalent to:
\begin{equation}
\label{e-l1-2}
\frac{1}{2}{\rm trace}({\cal B' \circ B}^{-1}(t))  \geq (n-1) \frac{s'_c(t)}{s_c(t)} 
\end{equation}
on all $(0,b_\xi )$. As ${\cal B}(t)$ is self-adjoint, we can choose an orthonormal basis $\{ e_1(t),...,e_{n-1}(t)\} $ of the (positive definite) subspace  $\xi^{\perp}$ consisting of eigenvectors of ${\cal B}(t)$. Calling $\lambda_1(t),...,\lambda_{n-1}(t)$ the corresponding eigenvalues, we have:
\begin{equation}
\label{e-l1-3}
{\rm trace}({\cal B' \circ B}^{-1}(t)) = \sum_{i=1}^{n-1} \frac{1}{\lambda_i(t)} g(e_i(t),{\cal B}'(t)e_i(t)) = \sum_{i=1}^{n-1} \frac{\frac{d}{ds}\mid_{s=t} \mid J^{\xi}_{e_i(t)}(s)\mid^{2}}
{\mid J^{\xi}_{e_i(t)}(t)\mid^{2}}. 
\end{equation}
Thus, to prove (\ref{e-l1-2}) it suffices to show the stronger inequality:
\begin{equation}
\label{e-l1-4}
\frac{\frac{d}{ds}\mid_{s=t} \mid J^{\xi}_{e_i(t)}(s)\mid^{2}}
{\mid J^{\xi}_{e_i(t)}(t)\mid^{2}} \geq 2 \frac{s'_c(t)}{s_c(t)}, 
\end{equation}
for each $i$. But, this inequality is just the consequence of applying the Lorentzian Rauch's comparison Theorem 
\ref{Rauch} to our case. The discussion of the equality is also straightforward from the equality in Theorem 
\ref{Rauch}. $\Box$

{\it Remark.} In above context in which a bound on sectional curvature is given, length estimates for Jacobi fields
are readily obtained by employing Rauch type comparison theorems.  Thus in Theorem \ref{t1} it is natural to 
consider the tensor field  ${\cal B}(t)  = \hat{{\cal A}}_{\xi}(t) \circ {\cal A}_{\xi}(t)$ as is commonly done in differential geometry nowadays, since
$g({\cal B}(t)(e_i),e_i) = g(J_{e_i}(t),J_{e_i}(t))$. On
the other hand, in obtaining the more delicate Bishop-Gromov type volume comparison results, only a Ricci 
curvature inequality is assumed and hence, the Jacobi equation is not so directly at hand.  It has been
accordingly helpful both in Riemannian geometry and in General Relativity to pass to associated matrix 
and scalar Ricatti equations by making the change of variables 
${\cal U}_{\xi}(t) = {\cal A}_{\xi}' \circ {\cal A}_{\xi}^{-1}(t)$
 as is done here in the 
proof of Lemma \ref{l2}. Inequality (\ref{Fi}) below may also be obtained as a by product of the more powerful 
Raychaudhuri equation standard in General Relativity (cf. \cite[Eq. (12.2), p. 430]{BEE}, for a textbook
exposition.)

\vspace{3mm}
\noindent {\bf Proof of Theorem \ref{t2}.} Reasoning as in the proof of Theorem \ref{t1}, and taking 
into account the Remark under Theorem \ref{t2}, the proof is reduced to the following result. 

\begin{lema}
\label{l2}
Put $\psi (t) \equiv {\det} {\cal A}_{\xi}(t) / s_c(t)^{n-1}$. Under the hypotheses of Theorem \ref{t2}:
\[
\psi(t)  \leq 1, \forall t \in (0,b_{\xi}).
\]
Moreover, if the equality holds at $t_0 \in (0,b_{\xi})$ then:
\[
 {\cal A}_{\xi}(t) = s_c(t) I, \quad {\rm and} \quad R_{\xi}(t)=-cI
\quad \forall t \in [0,t_0].
\]
\end{lema}

{\it Proof.} Again 
$\lim_{t\rightarrow 0} \psi (t)= 1$;  
thus, we will check just that $\psi'(t) \leq 0$ on $(0,b_{\xi})$, or, equivalently:
\begin{equation}
\label{e-l2-1}
\frac{d}{dt} {\rm \log (det}{\cal A}_{\xi})(t) \leq (n-1) \frac{s'_c(t)}{s_c(t)}. 
\end{equation}
The left member is equal to ${\rm trace}[{\cal A}_{\xi}'\circ {\cal A}_{\xi}^{-1}(t)]$. Consider as usual the 
following definitions:

\[
{\cal U}(t) \equiv {\cal U}_{\xi}(t) = {\cal A}_{\xi}' \circ {\cal A}_{\xi}^{-1}(t); \quad \Phi(t)= {\rm trace}\, {\cal U}\, (t)
\]
\[
\Phi_c(t) = (n-1) \frac{s'_c(t)}{s_c(t)} =(n-1) {\rm Ctg}_c (t)
\]
where ${\rm Ctg}_c (t) \equiv c_c (t) /s_c(t)$. Then, our objective (\ref{e-l2-1}) can be stated as:
\begin{equation}
\label{e-l2-2}
\Phi(t) \leq \Phi_c(t) \quad \forall t \in (0,b_{\xi}).
\end{equation}

\vspace{3mm}
\noindent As in the Riemannian case, ${\cal U}(t)$ satisfies the following properties:

(1) ${\cal U}(t)$ is self-adjoint, for all $t$.

(2) ${\cal U}(t)$ satisfies the matrix Riccati equation:
\begin{equation}
\label{e-l2-3}
{\cal U}' + {\cal U}^2 + {\cal R}_{\xi} = 0
\end{equation}

\vspace{3mm}
\noindent The problem is now reduced to scalar Riccati inequalities by using (1) and Cauchy-Schwarz inequality, 
$\Phi^2 \leq (n-1) {\rm trace}({\cal U}^2)$, and finally, taking the trace in (\ref{e-l2-3}):
\[
\Phi' + \frac{\Phi^2}{n-1} + {\rm Ric}(\gamma_{\xi}',\gamma_{\xi}') \leq 0.  
\]
As ${\rm Ric}(\gamma_{\xi}',\gamma_{\xi}') \geq (n-1) c \cdot 
g(\gamma_{\xi}',\gamma_{\xi}') = (-c)(n-1)$, we have putting $k=-c$:
\begin{equation}
\label{Fi}
\Phi' + \frac{\Phi^2}{n-1} + (n-1) k \leq 0.  
\end{equation}
On the other hand,  $\Phi$ satisfies the limit condition
\begin{equation}
\label{FiLi}
\lim_{t\rightarrow 0} t\Phi(t) = n-1
\end{equation}
(use that $\Phi(t) ={\rm det}{\cal A}'_{\xi}(t)/{\rm det}{\cal A}_{\xi}(t)$, and the property (2) of ${\cal A}_{\xi}$ in the proof of Theorem \ref{t1}). 

\vspace{3mm}
\noindent Finally, note that $\Phi_c$ is the solution of the Riccati equation:
\begin{equation}
\label{Fic}
\Phi'_c + \frac{\Phi_c^2}{n-1} + (n-1) k = 0  
\end{equation}
with limit condition 
\begin{equation}
\label{FiLic}
\lim_{t\rightarrow 0} t\Phi_c(t) = n-1.
\end{equation}
Thus, the standard comparison between the solutions to Riccati relations (\ref{Fi}), (\ref{Fic}) with equal limit conditions (\ref{FiLi}), (\ref{FiLic}) (see, for example, \cite[p. 121-2]{Ch}) yields the required inequality (\ref{e-l2-2}). $\Box$

\vspace{3mm}
\noindent {\bf Proof and discussion of Theorem \ref{t3}.} Clearly, the cut function of each $U^r$ satisfies 
$c_{U^r}=r c_U$ and:
\begin{equation}
\label{equation}
V(r)= \frac{\int_{\overline{U}[1]} d\xi \int_0^{rc_u(\xi)} dt \, \, {\rm det}({\cal A}_{\xi}(t))}{\int_{\overline{U}[1]} d\xi \int_0^{rc_u(\xi)} dt \, \, s_c^{n-1}(t)} = \frac{\int_{\overline{U}[1]} d\xi \int_0^{r} dt \, \, f(t,\xi)}{\int_{\overline{U}[1]} d\xi \int_0^{r} dt \, \, g(t,\xi)}
\end{equation}
where
\[
f(t,\xi)= c_U(\xi ) \cdot {\rm det}({\cal A}_{\xi}(c_U(\xi) t)) 
\]
\[
g(t,\xi)= c_U(\xi )  \cdot s_c^{n-1}(c_U(\xi) t)
\]      

\noindent Now, put
\[
F(t)=\int_{\overline{U}[1]} f(t,\xi ) \, \, d\xi 
\]
\[
G(t)=\int_{\overline{U}[1]} g(t,\xi ) \, \, d\xi .
\]

\noindent By \cite[Lemma 3.1, p. 124]{Ch}, 
{\it if the map $t\rightarrow F(t)/G(t)$ is non-increasing then so is $V(r)$}
and, thus Theorem \ref{t3} would be proven (note that the discussion of the equality would follow steps 
analogous to the Riemannian proof). This fact will be used to finish the proof of Cases {\it (A) {\rm and} (B)}.

\vspace{3mm}
\noindent {\it Proof of Case (A).} Using the expression for $s_c$ in the case that $c=0$: 
\[
G(t)= \int_{\overline{U}[1]} c_U(\xi)^n t^{n-1} d\xi = \lambda t^{n-1}
\]
with $\lambda = \int_{\overline{U}[1]} c_U(\xi)^n d\xi.$ Now, recall that, from the proof of Lemma \ref{l2}, $\psi (t) = {\rm det}({\cal A}_{\xi}(t))/t^{n-1} $ is non-increasing, and, thus, if $r<R, r,R\in I$:
\[
\lambda \frac{F(r)}{G(r)}= 
\int_{\overline{U}[1]} \frac{ c_u(\xi) \, \, {\rm det}({\cal A}_{\xi}(r))}{r^{n-1}} d\xi \geq 
\int_{\overline{U}[1]} \frac{ c_u(\xi) \, \, {\rm det}({\cal A}_{\xi}(R))}{R^{n-1}} d\xi = 
\lambda \frac{F(R)}{G(R)}
\]
as required.

\vspace{3mm}
\noindent {\it Proof of Case (B).} Assume $c_U(\xi)\equiv a \in \R{}$. Now 
\[
G(t)= a\lambda\cdot s_c^{n-1}(at)
\]
with $\lambda= \int_{\overline{U}[1]}  1 \cdot d\xi.$ Then, the result follows taking into account that $\psi(t) = {\rm det}({\cal A}_{\xi}(t))/s_c^{n-1} $ is non-increasing, as in Case (A). $\Box$

\vspace{3mm}
\noindent {\it Discussion for the necessity of either hypothesis (A) or (B).} If neither of the conditions  
{\it (A), (B)} are imposed, we have proven that, for each $\xi$, the quotient map 
\begin{equation}
\label{b1}
t\rightarrow f(t,\xi)/g(t,\xi)
\end{equation}
 is non-increasing and, thus, so is 
\begin{equation}
\label{b15}
r \rightarrow \int_0^r f(t,\xi)dt/\int_0^r g(t,\xi)dt.
\end{equation}
Nevertheless, we do not know if the function 
\begin{equation}
\label{b2}
t \rightarrow \int_{\overline{U}[1]} f(t,\xi)d\xi/\int_{\overline{U}[1]}  g(t,\xi)d\xi
\end{equation}
is non-increasing. In fact, it is not difficult to find functions $f,g$ satisfying the condition (\ref{b1}) is 
non-increasing, but (\ref{b2}) is not:

\begin{cont}
\label{3.3}
{\rm Taking into account the monotonicity of (\ref{b1}), our problem is equivalent to: given positive numbers $a_i, b_i, c_i, d_i, i\in\{ 1,...,k\}$ such that $a_i/b_i \geq c_i /d_i ,$ for all $i , $ is it true that 
\[ 
\sum_{i=1}^k a_i /\sum_{i=1}^k b_i \geq \sum_{i=1}^k c_i / \sum_{i=1}^k d_i \, ?
\]
[To check the equivalence, consider the integrals (\ref{b2}) as Riemann integrals, and take $k$ points $\{ \xi_1,...\xi_k \}$ in $\overline{U}[1]$ to make a Riemann integral sum; then put:
\[
a_i =f(r,\xi_i), \quad b_i=g(r,\xi_i), \quad c_i =f(R,\xi_i), \quad d_i=g(R,\xi_i).]
\]
Now, to answer negatively the question, just choose: 
\[
a_1 = a_2 = 2, b_1 = 1, b_2 = 1/M, c_1=1, c_2=M, d_1=d_2=1
\] 
for large $M$. Alternatively, putting:
\[
a_1 = 3/10, a_2 = 1/99, b_1 = 1/2, b_2 = 1/90, c_1=1/2, c_2=9/10, d_1=d_2=1
\] 
a counterexample is obtained satyisfying also: $a_i< c_i <d_i , a_i< b_i <d_i.$ So, we can identify:
\[
a_i = \int_0^r f(t,\xi_i) dt, \quad b_i= \int_0^r g(t,\xi) dt, \quad c_i =\int_0^R f(t,\xi_i) dt, \quad d_i=\int_0^R g(t,\xi) dt
\]
thus the relation corresponding to (\ref{b15}) is non-increasing, 
but nevertheless,  
\[
\frac{\int_{\overline{U}[1]} \int_0^r f }{\int_{\overline{U}[1]}  \int_0^r g} < \frac{\int_{\overline{U}[1]} \int_0^R f }{\int_{\overline{U}[1]}  \int_0^R g},
\]
that is, $V(r)<V(R)$. 
}
\end{cont}

\noindent Summing up, even though finding an explicit counterexample to the necessity of either $(A)$ or $(B)$ 
in Theorem \ref{t3} seems complicated, the previous computations suggest that this counterexample must exist. 
Note that it is suggested not only for the Lorentzian case but also if Theorem \ref{t3} is formulated for a 
SCV subset in the Riemannian case.

\section{Examples}

{\bf GRW spacetimes.} Consider a Generalized Robertson-Walker (GRW) spacetime  $(I\times F, g=-dt^2 + f^2(t) g_F)$, 
that is, a warped product with base $(I,-dt^2)$, $I\subseteq \R{}$ an open interval, the fiber any (connected) 
Riemannian manifold $(F,g_F)$ of dimension $m=n-1$, and warping function $f:I\rightarrow \R{}, f>0$. 
The elements of the fiber, like the Ricci or sectional curvature, will be denoted by adding a subscript 
$_F$ (${\rm Ric}_F , K_F$), and those of the GRW spacetime will have no index (Ric, $K$). In this Section 
we characterize when the curvatures of a GRW spacetime satisfy a bound which makes applicable the previous 
Theorems. Related properties about GRW space-times may be found in \cite{ARS-to}, \cite{grg3}, \cite{jmp3}.

\begin{prop}
\label{p2}
Fix $ p=(t,x)\in I\times F, c\in \R{}$. For a GRW spacetime, the condition:
\begin{equation}
\label{e-p2}
K(\pi) \geq c, \quad \hbox{ for all  timelike  planes }\, \pi \subset T_{p}(I\times F) 
\end{equation}
is equivalent to the following two conditions:

(A) $ f''(t)/f(t) \geq c$.

(B) $K_F (\pi_F ) \leq (f\cdot f''-f'^2)(t)$, for all tangent planes in  $T_xF$.

Moreover, the equality $K(\pi )=c$ holds for all timelike planes (and, then, for all nondegenerate planes) 
if and only if the warping function $f$ satisfies the equality in (A) and, when $m>1$, the curvature $K_F$ 
is a constant for which the equality in (B) holds.
\end{prop}

\begin{prop}
\label{p1}
Fix $p=(t,x)\in I\times F, c\in \R{}$. For a GRW spacetime, the condition:
\begin{equation}
\label{e-p1}
{\rm Ric}(X,X)\geq mc \cdot g(X,X), \quad \hbox{for all  timelike  vectors} \, X\in T_p (I\times F) 
\end{equation}
is equivalent to the following two conditions:

(A) $ f''(t)/f(t) \leq c$.

(B) ${\rm Ric}_F (Z,Z) \geq (m-1) (f\cdot f''-f'^2)(t)\cdot g_F(Z,Z)$, for all tangent vector $Z \in T_xF$.

Moreover, the equality holds in (\ref{e-p1}) for all timelike $X$ (and, then, the GRW spacetime is Einstein) 
if and only if the warping function $f$ satisfies the equality in (A) and, when $m>1$, the fiber is Einstein 
satisfying the equality in (B).

\end{prop}

In the case $m=1$ the conditions {\it (B)} are automatically satisfied. 
Of course, both propositions remain true if all the inequalities are reversed; as they are stated, examples for the Theorems in Section 2 are directly obtained. Clearly, (\ref{e-p2}) implies (\ref{e-p1}) reversed.
We will begin with the second result, which has a bit shorter proof.

\vspace{3mm}
\noindent {\it Proof of Proposition \ref{p1}}. Put 
$
X=\partial_t + \lambda Z
$
with $Z$ tangent to $F$ and $g_F$-unitary, $\lambda^2 f^2\leq 1$ (dependences on $p,t,x$ will be omitted). Then, 
clearly $g(X,X)=-1+\lambda^2f^2 (\leq 0)$ and, by  
\cite[p. 211]{O}: 
\begin{equation}
\label{ricci}
{\rm Ric}(X,X)= -m\frac{f''}{f} + \lambda^2\left( {\rm Ric}_F(Z,Z) + f\cdot f'' + (m-1)f'^2 \right) . 
\end{equation}
Thus, if (\ref{e-p1}) holds then inequalities {\it (A) {\em and} (B)} are obtained just by choosing $\lambda=0$ and (as a limit case) $\lambda = 1/f$. Conversely, if these two inequalies hold then, for each fixed $(t,x)\in I\times F$ and $g_F-$unitary $Z\in T_x F$,  define the functions 
\[
H_1(\lambda )= {\rm Ric} (\partial_t + \lambda Z,\partial_t + \lambda Z) = {\rm Ric}(\partial_t , \partial_t ) + \lambda^2 {\rm Ric}(Z,Z)
\]
\[
H_2(\lambda) = mc\cdot g(\partial_t + \lambda Z, \partial_t + \lambda Z) = -mc + \lambda^2 mc\cdot f^2
\]
\[
H(\lambda)=H_1(\lambda ) -H_2(\lambda ),
\] 
for all $\lambda \in [0,1/f]$; the proof ends by showing $H\geq 0$. Clearly, $H(\lambda ) \equiv a + b \lambda ^2 $ for certain $a,b \in \R{}$, that is, if non-constant then $H$ is a piece of parabola centered at the axis $\lambda \equiv 0$. But, by (\ref{ricci}) and hypotheses {\it (A), (B)}, $H$ satisfies at the extremes: $H(0)\geq 0, H(1/f)\geq 0$; thus $H\geq 0$ on all $[0,1/f]$, as required. $\Box$

\vspace{3mm}
\noindent {\it Proof of Proposition \ref{p2}}. Taking $Z$ unitary in $\pi \cap \partial_t^{\perp}$ a simple algebraic computation shows:

\begin{lema}
If $m\geq 2$ any timelike plane $\pi$ can be written as:
\[
\pi = {\rm Span}\{\partial_t + \lambda Y, Z\}
\]
where  $0\leq \lambda^2 f^2<1$ and $Y,Z$ are tangent to $F$ with $g(Y,Y)=g(Z,Z)=1, g(Y,Z)=0$.
\end{lema}

In what follows we will use this Lemma; the case $m=1$ is straightforward (and can be deduced from the next formula, just putting $\lambda=0, Y=Z$). By using \cite[p. 210]{O}: 
\[
K(\pi )= \frac{R(\partial_t , Z,Z, \partial_t ) + 2\lambda R(Y,Z,Z,\partial_t)+ \lambda ^2 R(Y,Z,Z,Y)}{(-1+\lambda^2 f^2)f^2}
\]
with:
\[
R(\partial_t , Z,Z, \partial_t ) =-f\cdot f'' 
\]
\[
R(Y,Z,Z,\partial_t) = 0 
\]
\[
R(Y,Z,Z,Y)= f^2\cdot K_F({\rm Span}\{Y,Z\}) + f^2\cdot f'^2.
\]
Thus, (\ref{e-p2}) is equivalent to:
\begin{equation}
\label{ep2}
-f\cdot f'' + \lambda^2 f^2 \left(K_F({\rm Span}\{Y,Z\}) + f'^2 \right) \leq -c f^2 + c \lambda^2 f^4 . 
\end{equation}
Putting $\lambda = 0$ or $\lambda = 1/f$ in (\ref{ep2}), the relations {\it (A), (B)} are obtained. For the converse, consider the functions $H_1(\lambda ) , H_2(\lambda )$ equal, respectively, to the left and right member of (\ref{ep2}). Then, it suffices to reason with the function $H(\lambda ) \equiv  H_1(\lambda )- H_2(\lambda )\, $ as in the previous  proof. $\Box$

\vspace{3mm}
\noindent {\bf Conformal changes.} Consider any Lorentzian manifold $(M,g)$, and a subset $U,$ SCLV at $p\in M$. Let $\omega : U\rightarrow \R{}$ be a {\it radial } function, that is, if $\xi, \xi' \in \overline{U}, r=$ $\mid\xi\mid =\mid \xi' \mid$ then $\omega (\exp_p (\xi )) = \omega (\exp_p (\xi' ))$. We will consider geodesic coordinates at $p$, and put $\omega \equiv \omega(r)$. Let us study the conformal metric
\[
g^* = e^{2\omega }\cdot g .
\]
Let $B=\{ v_1,v_2 \}$ be an orthonormal basis of a (non-degenerate) plane $\pi$ tangent to $U$, and put $\eta_i = g(v_i , v_i )$.  The relation between the sectional curvatures $K(\pi ), K^*(\pi )$ of $\pi$ for $g$ and $g^*$ is:

\begin{equation}
\label{conf}
e^{ 2 \omega } K^*(\pi) = K(\pi ) + \sum_{i=1}^2 \eta_i \left( v_i(\omega )^2 - {\rm Hess}\, \omega(v_i , v_i )\right) -g(\nabla \omega , \nabla \omega ) 
\end{equation}
(see, for example, \cite[p. 58]{Be}). If $v_1 = \partial/\partial r$ then: $v_1(\omega )= \omega', v_2(\omega )=0$, and $\sum_{i=1}^2 \eta_i  v_i(\omega )^2 -g(\nabla \omega , \nabla \omega ) = 0$. Thus, when $\pi$ is a (timelike) radial plane: 
\begin{equation}
\label{conf2}
e^{2\omega } K^*(\pi) = K(\pi )  + \omega'' - {\rm Hess}\, \omega(v_2 , v_2).
\end{equation}
If $\omega$ can be extended differentiably to 0, then, close to $p$, the term ${\rm Hess}\, \omega (v_2,v_2) (=-r\cdot \omega')$ is as close to 0 as we want. Thus, choosing, for example:
\begin{equation}
\label{confo}
\omega(r)= a\cdot r^2
\end{equation}
we have:
\begin{equation}
\label{conff}
\begin{array}{ll}
{\rm if}\quad a>0 \quad & K^*(\pi ) > K(\pi ) \\
{\rm if}\quad a<0 \quad & K^*(\pi ) < K(\pi )
\end{array}
\end{equation}
for all radial planes $\pi$ in a sufficiently small subset SCLV  at $p$. Inequalities (\ref{conff}) yield 
the required examples not only for Theorem \ref{t1} but also for Theorems \ref{t2}, \ref{t3}; nevertheless, 
for these last two results it is not difficult to make similar computations by considering the Ricci curvature
directly. In fact, note that
\[
{\rm Ric}^*(\partial_r, \partial_r ) = {\rm Ric}(\partial_r, \partial_r) - \Delta u - (n-2) u'', 
\]
where $\Delta$ denotes Laplacian and, close to 0, $\Delta u \cong u''$. 

\section{Local result for arbitrary semi-Riemannian manifolds}

First, our set-up will be extended to semi-Riemannian manifolds. Let $(M,g)$ be a semi-Riemannian 
$n-$manifold of arbitrary index, $n\geq 2$, and $Q_{c}$ the  model space of the same index and constant 
curvature $c$, for some fixed $c \in \R{}$;  choose $p\in M$ and $p_0 \in Q_c $, and choose a linear 
isometry $i: T_pM \rightarrow T_{p_0} Q_c $. Fix a normal neighborhood $V$ of $p$; 
$U$ will denote a subset which is SCV  at $p$ (as defined in Section 2), and we will always assume 
$U\subset V$. Now the transplantation $F$ is defined by using $V$, that is: $F:V\rightarrow Q_c, 
F=\exp_{p_0} \circ i \circ  (\exp_p \mid_V )^{-1}$. Put $\overline{V}_0 = i(\overline{V})$ and 
$V_0 = \exp_{p_0}(\overline{V}_0) (=F(V))$. $V$ will be chosen small enough such that $V_0$ is a 
normal neighborhood of $p_0$. Further, we will also compare $(M,g)$ with another semi-Riemannian 
manifold $(\hat{M},\hat{g})$ of equal index and dimension. We will also fix $\hat{p} \in \hat{M}$ 
and a linear isometry $\hat{i}:T_pM\rightarrow T_{\hat{p}}\hat{M}$; the symbol $\hat{}$ will be 
put on the elements of $\hat{M}$ or on those elements necessary for the corresponding comparison, like 
the Ricci tensor $\hat{{\rm Ric}}$ or 
$\hat{F}=\exp_{\hat{p}} \circ \hat{i} \circ  (\exp_p \mid_V )^{-1}, \hat{V}=\hat{F}(V)$.

\begin{lema}
\label{l}
Assume that 
\begin{equation}
\label{e-l}
{\rm Ric}(\xi,\xi) > \hat{{\rm Ric}}(\hat{\xi},\hat{\xi}) \quad ({\rm resp.} \quad {\rm Ric}(\xi,\xi) < \hat{{\rm Ric}}(\hat{\xi},\hat{\xi}) ) 
\end{equation}
for some unit $\xi \in T_pM$, $\hat{\xi} = d\hat{F}_p (\xi) \in T_{\hat{p}}\hat{M}. $ 
Then there exists an open subset $U,$ SCV at $p,$ such that 
\[
{\rm vol}(A) < {\rm vol}(\hat{A} ) \quad ({\rm resp.} \quad {\rm vol}(A) > {\rm vol}(\hat{A} ))
\]
for all measurable subsets $A\subseteq U, \hat{A}=F(A)$.
\end{lema}
We postpone the proof of this result, 
and explore first its consequences.

\begin{teor}
\label{t}
Assume that inequalities {\rm (\ref{e-l})} holds. 
Then there exist an arbitrarily small normal neighborhood $W$ of $p$ such that 
\[
{\rm vol}(W) < {\rm vol}(\hat{W} ) \quad ({\rm resp.} \quad {\rm vol}(W) > {\rm vol}(\hat{W} )). 
\]
\end{teor}
[We mean by ``arbitrarily small" (for a class of neighborhoods at $p$) the possibility of choosing a topological basis at $p$ of such neighborhoods.]

{\it Proof.} Applying Lemma \ref{l} with $A=U$, the strict inequality in the volumes allows to enlarge $U$ slightly to obtain a normal neighborhood where the inequality between the volumes is preserved. $\Box$

For metrics of constant curvature, choosing a spacelike $\xi$ as well as a timelike one, both equalities in  Theorem \ref{t} are obtained. More precisely:

\begin{coro}
\label{c}
Take $c, \hat{c}\in \R{}, c\neq \hat{c}$ and the corresponding indefinite model spaces $Q_c , Q_{\hat{c}}$. Then, for each $p\in Q_c$, $\hat{p} \in Q_{\hat{c}}$ there are two arbitrarily small normal neighborhoods $W,W'$ of $p$ such that 
\begin{equation}
\label{*}
\begin{array}{c}
{\rm vol}(W) < {\rm vol}(\hat{W} ), \\ 
{\rm vol}(W') > {\rm vol}(\hat{W}' ), 
\end{array}
\end{equation}
where $\hat{W} = \hat{F}(W), \hat{W}' = \hat{F}(W')$.
\end{coro}

\vspace{3mm}

\noindent {\it Remark.} (1) When the metric $g$ is reversed to $-g$ then the sectional curvature changes sign. This makes consistent the comparison between the neighborhoods obtained for $Q_{\hat{c}}$ and $Q_c$ (using a $\xi$ either spacelike or timelike) with the comparison (reversing metrics) of $Q_{-\hat{c}}$ and $Q_{-c}$ (taking $\xi$ either timelike or spacelike).

(2) As a consequence of Corollary \ref{c}, the conclusion of the Theorem \ref{t} {\it cannot} be strengthened 
in the sense that the conclusion of Lemma \ref{l} holds for more arbitrary neighborhoods of p. That is, it is 
{\it not} true for any (normal) neighborhood $W$ of $p$ that ${\rm vol}(A) < {\rm vol}(\hat{A} )$ for all 
measurable subsets $A \subseteq W$. 

(3) Clearly, the conclusion of Corollary \ref{c} also holds replacing the model spaces by Lorentzian manifolds of dimension 2, and regarding now  $c, \hat{c}$ as their curvatures at $p$. 

\vspace{3mm}

\noindent  As a consequence of Lemma \ref{l} and Theorem \ref{t}, local versions of  Theorems \ref{t1}, and 
\ref{t2} are obtained. Note that in the local version of Theorem \ref{t1} the assumption on the sectional 
curvature is replaced by a condition just on the Ricci curvature. This remains true in the Riemannian case, 
and the following local version of Günther's theorem can be obtained, where the assumption on the sectional 
curvature is replaced by one on the Ricci curvature at a point (the proof is done at the end of the section).

\begin{coro}
\label{c-tf}
Assume that $(M,g), (\hat{M},\hat{g})$ are Riemannian. If 
\begin{equation}
\label{e-GB}
{\rm Ric}(\xi,\xi) < \hat{{\rm Ric}}(\hat{\xi},\hat{\xi}) 
\end{equation}
for all unit $\xi \in T_pM$, $\hat{\xi} = d\hat{F}_p (\xi) \in T_{\hat{p}}\hat{M}, $ 
then there exist a normal neighborhood $W$ of $p$ such that 
\[
{\rm vol}(A) > {\rm vol}(\hat{A} )
\]
for all measurable subsets $A\subseteq W, \hat{A}=F(A)$. 

In particular, there exist $\epsilon >0$ such that for all $r\in (0,\epsilon)$ the volume of the metric ball of radius $r$ centered at $p$, $B(p,r),$ is strictly bigger than the volume of the corresponding metric ball  at $\hat{p}$, $B(\hat{p},r)$.
\end{coro}

\vspace{3mm}
\noindent {\bf Proof of Lemma \ref{l}}. Consider an orthonormal basis $\{ e_1,\dots , e_{n-1}\}$ of $\xi^{\perp}$ and put $\pi_i = {\rm Span}\{\xi , e_i\}$, with sectional curvature $K(\pi_i)$. Let $\{ E_1,\dots , E_{n-1} \} $ be an orthonormal parallel basis of $\gamma_{\xi}^{'\perp}$ with $E_i(0)=e_i$ and $ \{ J_1,\dots , J_{n-1} \}$ Jacobi fields along $\gamma_{\xi}$ with $J_i(0)=0, J'_i(0)=e_i$,  for all $i$. Put $m=n-1$, $\epsilon = g(\xi, \xi) (=\pm 1)$. 
and $\epsilon_i = g(E_i,E_i)$, let $\delta_{ij}$ be Kronecker's delta and denote by ${\cal O}(t^s)$ terms of order $s$ in the corresponding series expansion.

\begin{lema}
\label{lA}
$J_i(t)$ admits the expansion:
\[
J_i(t) = \left(t- \frac{\epsilon}{6} K(\pi_i) \, t^3 + {\cal O}_i (t^4)\right) E_i(t) + \sum_{j\neq i}^{m} {\cal O}_{ji}(t^3)\, E_j(t).
\]
\end{lema}
{\it Proof.} Clearly, 
\[
\begin{array}{c}
g(J_i,E_j)(0)=0, \quad \quad g(J_i,E_j)'(0) \, (= g(J'_i,E_j)(0)) = \epsilon_i \delta_{ij} \\

g(J_i,E_j)''(0)= -R(\gamma'_{\xi},J_i,E_j,\gamma'_{\xi}))(0)=0 \quad \quad  g(J_i,E_j)'''(0)= -R(\gamma'_{\xi},J'_i,E_j,\gamma'_{\xi}))(0). 
\end{array}
\]

\noindent Thus, for $i=j$,  $g(J_i,E_i)'''(0)= -\epsilon \, \epsilon_i\, K(\pi_i)$, and the result follows using Taylor's expansion in 
$J_i (t) = \sum_{j=1}^{m} \epsilon_j g(J_i,E_j)(t)\, E_j(t).$ $\Box$ 

Consider the endomorphism field ${\cal A}_{\xi}$ along $\gamma_{\xi}$ in Section 3. Lemma \ref{lA} will be 
used to prove the following result.

\begin{lema}
\label{lB}
${\rm det}{\cal A}_{\xi}(t)$ admits the expansion:
\[
{\rm det}{\cal A}_{\xi}(t) = t^m -\frac{1}{6} {\rm Ric}(\xi,\xi) \, t^{m+2} + {\cal O} (t^{m+3}).
\]
\end{lema}
{\it Proof.} For $m=1$ the result is trivial; otherwise, 
${\rm det}{\cal A}_{\xi}(t) $
\[ 
= {\rm det}\left( \begin{array}{cccc}
t - \frac{\epsilon}{6} K(\pi_1) \, t^3 + {\cal O}_1 (t^4)  &  
{\cal O}_{12}(t^3) &\cdots & {\cal O}_{1m}(t^3) \\
{\cal O}_{21}(t^3) & 
t- \frac{\epsilon}{6} K(\pi_2) \, t^3 + {\cal O}_2 (t^4)  &  
\cdots & {\cal O}_{2m}(t^3) \\
\vdots & \vdots & & \vdots \\
{\cal O}_{m1}(t^3) & {\cal O}_{m2}(t^3) & \cdots & 
t- \frac{\epsilon}{6} K(\pi_m) \, t^3 + {\cal O}_m (t^4)
\end{array} \right) 
\]

\[
= \prod_{i=1}^{m} \left( t- \frac{\epsilon}{6} K(\pi_i) \, t^3\right) + {\hat {\cal O}}(t^{m+3}), 
\]
and taking into account ${\rm Ric}(\xi, \xi) = \epsilon \sum_{i=1}^m K(\pi_i)$ the result follows. $\Box$

Now, Lemma \ref{l} can be easily proven as follows. Consider the endomorphism field ${\cal A}_{\hat{\xi}}$ along $\gamma_{\hat{\xi}}, $ and $\hat{\xi} = \hat{i}(\xi )$. By Lemma \ref{lB},  if $ {\rm Ric}(\xi, \xi) > {\rm \hat{Ric}}(\hat{\xi}, \hat{\xi})$ then there exists a $\delta >0$ such that:
\[
{\rm det}{\cal A}_{\xi}(t) < {\rm det}{\cal A}_{\hat{\xi}}(t), \quad \forall t\in (0,\delta).
\]
By continuity, this inequality holds for radial geodesics close to $\xi$ on a subset $U$ SCV  at $p$, and an 
integration in normal coordinates yields the result. 

\vspace{3mm}

\noindent {\bf Proof of Corollary \ref{c-tf}}. Consider any normal neighborhood $V$ of $p$, and define the function $c: S_pM \rightarrow \R{}, c(\xi )= {\rm Sup}\{t\in \R{} \mid t\xi \in V, \quad {\rm and} \quad {\rm det} {\cal A}_{\xi}(t) \geq {\rm det}{\cal A}_{\hat{\xi}}(t)\}$, where $S_pM $ is the unit sphere at $T_pM$. The result follows putting $\epsilon$ equal to  the minimum of this function (which, by Lemma \ref{lB}, is strictly positive) and:
\[
W=\{ \exp_p(t\xi) \mid \xi \in S_pM, 0\leq t < c(\xi )\}. \quad \Box
\]

\vspace*{2mm}

\begin{tabbing}
Department of Mathematics \\
358 Little Hall, \\
PO Box 11,8105, \\
University of Florida \\
Gainesville FL 32611-8105, U.S.A. \\
\end{tabbing}

\begin{tabbing}
Departamento de Geometría y Topología \\
Fac. Ciencias, Avda. Fuentenueva s/n, \\
Universidad de Granada, \\
18071 Granada, Spain. \\
\end{tabbing}

\newpage

\vspace*{16cm}
\begin{center}
 Figure 1\\
Transplantation $F$ of a subset SCLV at $p$.\\
\end{center}

\end{document}